\documentstyle[12pt]{article}
\newcommand{\too}{\longrightarrow}
\newcommand{\om}{\omega}

\newcommand{\T}{{\cal T}}

\newcommand{\s}{{\cal S}}

\newcommand{\p}{{\cal P}}

\newcommand{\di}{\displaystyle}
\newcommand{\Om}{\Omega}
\newcommand{\na}{\nabla}

\newcommand{\al}{\alpha}
\newcommand{\be}{\beta}

\newcommand{\Ga}{\Gamma}

\newcommand{\la}{\lambda}
\newcommand{\De}{\Delta}
\newcommand{\de}{\delta}

\def \reel{ {\rm I}\!{\rm R} }

 \def \comp{ \;{}^{ {}_\vert }\!\!\!{\rm C} }

 \def \nat{ { {\rm I}\!{\rm N}} }

 \def \rat{ {\rm Q}\kern-.65em {}^{{}_/ }}

\newtheorem{th}{Theorem}[section]

\newtheorem{Le}{Lemma}[section]
\newtheorem{co}{Corollary}[section]
\newtheorem{rem}{Remark}[section]
\title{ Spectra and symmetric eigentensors of the Lichnerowicz Laplacian
 on $S^n$ }
\author{M. Boucetta\\
Facult\'e des Sciences et Techniques \\
BP 549 Marrakech\\
Morocco
\\
Email: {\it boucetta@fstg-marrakech.ac.ma } \footnote{Recherche
men\'ee dans le cadre du Programme
 Th\'ematique d'Appui \`a la Recherche Scientifique PROTARS III.}} \date{ }\parindent=0cm
\begin{document}
\maketitle

{\bf Abstract.} We compute the eigenvalues with multiplicities  of
the Lichnerowicz Laplacian acting on the space of symmetric
covariant tensor fields on the Euclidian sphere $S^n$. The spaces
of symmetric eigentensors are explicitly given.

\bigskip



\section{Introduction}
Let $(M,g)$ be a Riemannian $n$-manifold. For any $p\in\nat$, we
shall denote by $\Ga(\otimes^pT^*M)$, $\Om^p(M)$ and $\s^pM$ the
space of covariant $p$-tensor fields on $M$, the space of
differential $p$-forms  on $M$ and the space of symmetric
covariant $p$-tensor fields on $M$, respectively. Note that
$\Ga(\otimes^0T^*M)=\Om^0(M)=\s^0M=C^\infty(M,\reel)$,
$\di\Om(M)=\sum_{p=0}^n\Om^p(M)$ and
$\di\s(M)=\sum_{p\geq0}\s^p(M)$.

Let $D$ be the Levi-Civita connection associated to $g$;  its
curvature tensor field $R$ is given by
\[R(X,Y)Z=D_{[X,Y]}Z-\left(D_XD_YZ-D_YD_XZ\right),\]and
the Ricci endomorphism field $r:TM\too TM$
  is given by
\[g(r(X),Y)=\sum_{i=1}^ng(R(X,E_i)Y,E_i),\]where
$(E_1,\ldots,E_n)$ is any local orthonormal frame.

 For any $p\in\nat$, the connection $D$ induces a
differential operator
$D:\Ga(\otimes^pT^*M)\too\Ga(\otimes^{p+1}T^*M)$ given by
\[DT(X,Y_1,\ldots,Y_p)=D_XT(Y_1,\ldots,Y_p)=X.T(Y_1,\ldots,Y_p)-
\sum_{j=1}^pT(Y_1,\ldots,D_XY_j,\ldots,Y_p).\]  Its formal adjoint
$D^*:\Ga(\otimes^{p+1}T^*M)\too\Ga(\otimes^{p}T^*M)$ is  given by
\[D^*T(Y_1,\ldots,Y_p)=-\sum_{j=1}^nD_{E_i}T(E_i,Y_1,\ldots,Y_p),\]
where $(E_1,\ldots,E_n)$ is any local orthonormal frame.

Recall that, for any differential $p$-form $\al$, we have
\begin{equation}
d\al(X_1,\ldots,X_{p+1})=\sum_{j=1}^{p+1}(-1)^{j+1}D_{X_j}\al(X_1,\ldots,
\hat{X}_j,\ldots,X_{p+1}).\end{equation}

We  denote by $\de$  the restriction of $D^*$ to
$\Om(M)\oplus\s(M)$ and we define $\de^*:\s^p(M)\too\s^{p+1}(M)$
by
$$\de^*T(X_1,\ldots,X_{p+1})=\sum_{j=1}^{p+1}D_{X_j}T(X_1,\ldots,
\hat{X}_j,\ldots,X_{p+1}).$$

Recall that the operator trace $Tr:\s^p(M)\too\s^{p-2}(M)$ is
given by
$$TrT(X_1,\ldots,X_{p-2})=\sum_{j=1}^nT(E_j,E_j,X_1,\ldots,X_{p-2}),$$
where $(E_1,\ldots,E_n)$ is any local orthonormal frame.

The Lichnerowicz Laplacian is the second order differential
operator $$\De_M:\Ga(\otimes^{p}T^*M)\too\Ga(\otimes^{p}T^*M)$$
given by
\[\De_M(T)=D^*D(T)+R(T),\] where $R(T)$ is the curvature
operator given by
\begin{eqnarray*}
&&R(T)(Y_1,\ldots,Y_p)=\sum_{j=1}^pT
(Y_1,\ldots,r(Y_j),\ldots,Y_p)\\
&-&\sum_{i<j}\sum_{l=1}^n\left\{T(Y_1,\ldots,E_l,\ldots,R(Y_i,E_l)Y_j,\ldots,Y_p)
+
T(Y_1,\ldots,R(Y_j,E_l)Y_i,\ldots,E_l,\ldots,Y_p)\right\},\nonumber\end{eqnarray*}
where $(E_1,\ldots,E_n)$ is any local orthonormal frame and, in
$$T(Y_1,\ldots,E_l,\ldots,R(Y_i,E_l)Y_j,\ldots,Y_p),$$ $E_l$ takes
the place of $Y_i$ and $R(Y_i,E_l)Y_j$ takes the place of $Y_j$.

This differential operator, introduced by Lichnerowicz in [15] pp.
26, is self-adjoint, elliptic and respects the symmetries of
tensor fields. In particular, $\De_M$ leaves invariant $\s(M)$
 and the restriction of
$\De_M$ to  $\Om(M)$ coincides with the Hodge-de Rham Laplacian,
i.e., for any differential $p$-form $\al$,
\begin{equation}\De_M\al=(d\de+\de d)(\al).\end{equation}

We have shown in [6] that, for any symmetric covariant tensor
field $T$,
\begin{equation}
\De_M(T)=(\de\circ\de^*-\de^*\de)(T)+2R(T).\end{equation}

Note that if $T\in\s(M)$ and $g^l$ denotes the symmetric product
of $l$ copies of the Riemannian metric $g$, we have
\begin{eqnarray}(Tr\circ\De_M) T&=&(\De_M\circ Tr) T,\\
\De_{M}(T\odot g^l)&=&(\De_MT)\odot g^l,\end{eqnarray}where
$\odot$ is the symmetric product.

The Lichnerowicz Laplacian  acting on symmetric covariant tensor
fields is of fundamental importance in mathematical physics  (see
for instance  [9], [20] and [22]). Note also that the Lichenrowicz
Laplacian acting on symmetric covariant 2-tensor fields appears in
many problems in Riemannian geometry (see [3], [5], [19]...).

On a Riemannian compact manifold, the Lichnerowicz Laplacian
$\De_M$ has discrete eigenvalues with finite
 multiplicities.
 For a given Riemannian compact manifold, it may be an interesting problem
 to determine explicitly
the eigenvalues and the eigentensors of $\De_M$ on $M$.

Let us enumerate the cases where the spectra of $\De_M$ was
computed:
\begin{enumerate}
\item $\De_M$ acting on $C^\infty(M,\comp)$:  $M$ is either flat
toris or
 Klein bottles [4], $M$ is a Hopf manifolds
 [1];
\item $\De_M$ acting on  $\Om(M)$: $M=S^n$ or $P^n(\comp)$ [10]
and [11], $M=\comp aP^2$ or $G_2/SO(4)$ [16] and [18],
$M=SO(n+1)/SO(2)\times SO(n)$ or $M=Sp(n+1)/Sp(1)\times Sp(n)$
[21]; \item $\De_M$ acting on $\s^2(M)$ and $M$ is the complex
projective space $P^2(\comp)$ [22]; \item $\De_M$ acting on
$\s^2(M)$ and $M$ is either $S^n$ or
 $P^n(\comp)$ [6] and [7];
 \item  Brian and Richard Millman give in [2] a theoretical
 method for
computing the spectra of Lichnerowicz Laplacian acting on $\Om(G)$
where $G$ is a compact semisimple Lie group endowed with the
Killing form; \item Some partial results where given in
[12]-[14].\end{enumerate}

In this paper, we compute the eigenvalues and we determine the
spaces of eigentensors of $\De_M$ acting on $\s(M)$ in the case
where $M$ is the Euclidian sphere $S^n$.

 Let us describe  our method briefly. We consider the
 $(n+1)$-Euclidian space $\reel^{n+1}$ with its canonical coordinates $(x_1,\ldots,x_{n+1})$.
  For any $k,p\in\nat$, we denote by
$\s^pH_k^\de$ the space of symmetric covariant $p$-tensor fields
$T$ on $\reel^{n+1}$ satisfying:
\begin{enumerate}\item $\di T=\sum_{1\leq i_1\leq\ldots\leq i_p\leq
n+1}T_{i_1,\ldots,i_p}dx_{i_1}\odot\ldots\odot dx_{i_p}$ where
$T_{i_1,\ldots,i_p}$ are homogeneous polynomials of degree $k$;
\item $\de(T)=\De_{\reel^{n+1}}(T)=0.$\end{enumerate}

 The $n$-dimensional sphere $S^n$ is the space of unitary
vectors in  $\reel^{n+1}$ and  the Euclidian metric on
$\reel^{n+1}$ induces a Riemannian metric on $S^n$. We denote by
$i:S^n\hookrightarrow\reel^{n+1}$ the canonical inclusion.

For any tensor field $T\in\Ga(\otimes^pT^*\reel^{n+1})$, we
compute $i^*(\De_{\reel^{n+1}}T)-\De_{S^n}(i^*T)$  and get a
formula (see Theorem 2.1). Inspired by this formula and having in
mind  the fact that $i^*:\di\sum_{k\geq0}\s^pH_k^\de\too\s^pS^n$
is injective and its image is dense in $\s^pS^n$ (see [10]), we
give, for any $k$, a direct sum decomposition of $\s^pH_k^\de$
composed by eigenspaces of $\De_{S^n}$. Thus we obtain the
eigenvalues and the spaces of eigentensors with its multiplicities
of $\De_{S^n}$ acting on $\s(S^n)$ (see Section 4).

Note that the eigenvalues and the eigenspaces of $\De_{S^n}$
acting on $\Om(S^n)$ was computed in [10] by using the
representation theory. In [11], I. Iwasaki and K. Katase recover
the result by a method  using the restriction of harmonic tensor
fields and a result in [8]. The formula obtained in Theorem 2.1
combined with the methods developed in [10] and [11] permit to
present those results in a more precise form (see Section 3).

\section{A relation between $\De_{\reel^{n+1}}$ and $\De_{S^n}$}

We consider the Euclidian space $\reel^{n+1}$ endowed with its
canonical coordinates $(x_1,\ldots,x_{n+1})$ and its canonical
Euclidian flat Riemannian metric $<\;,\;>$.  We denote by $D$ be
the Levi-Civita covariant derivative  associated to $<\;,\;>$. We
consider the radial vector field given by
$$\overrightarrow{r}=\sum_{i=1}^{n+1}x_i\frac{\partial}{\partial
x_i}.$$

 For any $p$-tensor field $T\in\Ga(\otimes^pT^*\reel^{n+1})$ and for any $1\leq
 i<j\leq p$, we denote by $i_{\overrightarrow{r},j}T$ the
$(p-1)$-tensor field given by
\[
i_{\overrightarrow{r},j}T(X_1,\ldots,X_{p-1})=T(X_1,\ldots,X_{j-1},
\overrightarrow{r},X_j,\ldots, X_{p-1}),
\]and by $Tr_{i,j}T$ the $(p-2)$-tensor field given by
\[ Tr_{i,j}T(X_1,\ldots,X_{p-2})=\sum_{l=1}^{n+1}T(X_1,\ldots,X_{i-1},
E_l,X_i,\ldots,X_{j-2},E_l,X_{j-1},\ldots, X_{p-2}),\]where
$(E_1,\ldots,E_{n+1})$ is any orthonormal basis of $\reel^{n+1}$.
Note that $Tr_{i,j}T=0$ if $T$ is a differential form and
$Tr_{i,j}T=TrT$ if $T$ is symmetric.

For any permutation $\sigma$ of $\{1,\ldots,p\}$, we denote by
$T^\sigma$ the $p$-tensor field
\[
T^\sigma(X_1,\ldots,X_{p})=T(X_{\sigma(1)},\ldots,X_{\sigma(p)}).\]
For $1\leq i<j\leq p$, the transposition of $(i,j)$ is the
permutation $\sigma_{i,j}$ of $\{1,\ldots,p\}$ such that
$\sigma_{i,j}(i)=j$, $\sigma_{i,j}(j)=i$ and $\sigma_{i,j}(k)=k$
for $k\not=i,j$. Let $\T$ denote the set of the transpositions of
$\{1,\ldots,p\}$.

The sphere $i:S^n\hookrightarrow\reel^{n+1}$ is endowed with the
Euclidian metric.

\begin{th} Let $T$ be a covariant $p$-tensor field on $\reel^{n+1}$. Then,
\begin{eqnarray*}
i^*(\De_{\reel^{n+1}}T)&=&\De_{S^n}i^*T+i^*\left(p(1-p)T+(2p-n+1)L_{\overrightarrow{r}}T-
L_{\overrightarrow{r}}\circ
L_{\overrightarrow{r}}T\right.\\
&&\left.-2\sum_{\sigma\in\T}T^\sigma+O(T)\right),\end{eqnarray*}
where $O(T)$ is given by
\begin{eqnarray*}
O(T)(X_1,\ldots,X_p)&=&2\sum_{i<j}<X_i,X_j>Tr_{i,j}(X_1,\ldots,\hat{X}_i,\ldots,\hat{X}_j,
\ldots,X_p)\\
&&-2\sum_{j=1}^pD_{X_j}(i_{\overrightarrow{r},j}T)(X_1,\ldots,\hat{X_j},\ldots,X_p),
\end{eqnarray*}where $\hat{X}$ designs that $X$ is deleted.
\end{th}

{\bf Proof.} The proof is a massive computation in a local
orthonormal frame using the properties of the Riemannian embedding
of the sphere in the Euclidian space.

We choose a local orthonormal frame of $\reel^{n+1}$
 of the form $(E_1,\ldots,E_n,N)$ such
that $E_i$ is tangent to $S^n$ for $1\leq i\leq n$ and
$N=\frac1r\overrightarrow{r}$ where
$r=\sqrt{x_1^2+\ldots+x_{n+1}^2}$.

For any vector field $X$ on $\reel^{n+1}$, we have
\begin{eqnarray} D_XN&=&\frac1r\left(X-<X,N>N\right),\\
D_NX&=&[N,X]+\frac1r(X-<X,N>N).\end{eqnarray}

Let $\na$ be the Levi-Civita connexion of the Riemannian metric on
$S^n$. We have, for any vector fields  $X,Y$ tangent to $S^n$,
\begin{equation} D_XY=\na_XY-<X,Y>N.\end{equation}

Let $T$ be a covariant $p$-tensor field on $\reel^{n+1}$ and
$(X_1,\ldots,X_p)$ a family of vector fields on $\reel^{n+1}$
which are tangent to $S^n$. A direct calculation using the
definition of the Lichnerowicz Laplacian gives$$
\begin{array}{l}
\di\De_{\reel^{n+1}}(T)(X_1,\ldots,X_p)=
D^*D(T)(X_1,\ldots,X_p)\\\di=\sum_{i=1}^{n}\left(-E_iE_i.T(X_1,\ldots,X_p)
+2\sum_{j=1}^pE_i.T(X_1,\ldots,D_{E_i}X_j,\ldots,X_p)\right.\\
\di+D_{E_i}E_i.T(X_1,\ldots,X_p)
-\sum_{j=1}^pT(X_1,\ldots,D_{D_{E_i}E_i}X_j,\ldots,X_p)\\
\di-\sum_{j=1}^pT(X_1,\ldots,D_{E_i}D_{E_i}X_j,\ldots,X_p)
\di\left.-2\sum_{l<j}T(X_1,\ldots,D_{E_i}X_l,\ldots,D_{E_i}X_j,\ldots,X_p)\right)
\end{array}$$
$$\begin{array}{l} -N.N.T(X_1,\ldots,X_p)
\di+2\sum_{j=1}^pN.T(X_1,\ldots,D_{N}X_j,\ldots,X_p)\\
+D_{N}N.T(X_1,\ldots,X_p)
\di-\sum_{j=1}^pT(X_1,\ldots,D_{D_{N}N}X_j,\ldots,X_p)\\
\di-\sum_{j=1}^pT(X_1,\ldots,D_{N}D_{N}X_j,\ldots,X_p)
-2\sum_{l<j}T(X_1,\ldots,D_{N}X_l,\ldots,D_{N}X_j,\ldots,X_p).\end{array}$$

 (6)-(8) make it obvious that
\begin{eqnarray}
D_{D_{E_i}E_i}X_j&=&\na_{\na_{E_i}E_i}X_j-<\na_{E_i}E_i,X_j>N-[N,X_j]\\
&&-\frac1r(X_j-<X_j,N>N),\nonumber\\
D_{E_i}D_{E_i}X_j&=&\na_{E_i}\na_{E_i}X_j-(<E_i,\na_{E_i}X_j>+E_i.<E_i,X_j>)N
\nonumber\\&&-
\frac1r<E_i,X_j>E_i,\\
D_ND_NX&=&[N,[N,X]]+\frac2r[N,X]+(\frac1{r^2}-\frac1r)(X-<X,N>N)\nonumber\\&&-\frac2rN.<X,N>N.
\end{eqnarray}

By  (8)-(10), we get easily, in restriction to $S^n$,
$$
\begin{array}{l}
\di\sum_{i=1}^n\left(2\sum_{j=1}^pE_i.T(X_1,\ldots,D_{E_i}X_j,\ldots,X_p)
\di+D_{E_i}E_i.T(X_1,\ldots,X_p)\right.\\\di
\left.-\sum_{j=1}^pT(X_1,\ldots,D_{D_{E_i}E_i}X_j,\ldots,X_p)
\di-\sum_{j=1}^pT(X_1,\ldots,D_{E_i}D_{E_i}X_j,\ldots,X_p)\right)=\\\di
\di\sum_{i=1}^n\left(2\sum_{j=1}^pE_i.T(X_1,\ldots,\na_{E_i}X_j,\ldots,X_p)
\di+\na_{E_i}E_i.T(X_1,\ldots,X_p)\right.\\\di
\left.-\sum_{j=1}^pT(X_1,\ldots,\na_{\na_{E_i}E_i}X_j,\ldots,X_p)
\di-\sum_{j=1}^pT(X_1,\ldots,\na_{E_i}\na_{E_i}X_j,\ldots,X_p)\right)\\\di
-2\sum_{j=1}^pX_j.T(X_1,\ldots,\stackrel{j}{\overbrace{N}},\ldots,X_p)+
p(n+1)T(X_1,\ldots,X_p)-nL_NT(X_1,\ldots,X_p).
\end{array}$$
On other hand, also by using (8), we have
$$\begin{array}{l}
\di\sum_{l<j}\sum_{i=1}^n
T(X_1,\ldots,D_{E_i}X_l,\ldots,D_{E_i}X_j,\ldots,X_p)=\\\di
\sum_{l<j}\sum_{i=1}^n
T(X_1,\ldots,D_{E_i}X_l,\ldots,\na_{E_i}X_j,\ldots,X_p)
-\sum_{l<j}T(X_1,\ldots,D_{X_j}X_l,\ldots,\stackrel{j}{\overbrace{N}},\ldots,X_p)=\\\di
\sum_{l<j}\sum_{i=1}^n
T(X_1,\ldots,\na_{E_i}X_l,\ldots,\na_{E_i}X_j,\ldots,X_p)-
\sum_{l<j}
T(X_1,\ldots,\stackrel{l}{\overbrace{N}},\ldots,\na_{X_l}X_j,\ldots,X_p)\\\di
-\sum_{l<j}T(X_1,\ldots,D_{X_j}X_l,\ldots,\stackrel{j}{\overbrace{N}},\ldots,X_p)=
\sum_{l<j}\sum_{i=1}^n
T(X_1,\ldots,\na_{E_i}X_l,\ldots,\na_{E_i}X_j,\ldots,X_p)\\\di
-\sum_{l<j}T(X_1,\ldots,D_{X_j}X_l,\ldots,\stackrel{j}{\overbrace{N}},\ldots,X_p)
-\sum_{l<j}
T(X_1,\ldots,\stackrel{l}{\overbrace{N}},\ldots,D_{X_l}X_j,\ldots,X_p)\\\di
-\sum_{l<j}<X_l,X_j>T(X_1,\ldots,\stackrel{l}{\overbrace{N}},\ldots,
\stackrel{j}{\overbrace{N}},\ldots,X_p).
\end{array}$$

So we get, in restriction to $S^n$, since $D_NN=0$
$$\begin{array}{l}\di\De_{\reel^{n+1}}(X_1,\ldots,X_p)-\na^*\na
T(X_1,\ldots,X_p)=\\\di
p(n+1)T(X_1,\ldots,X_p)-nL_NT(X_1,\ldots,X_p)-2\sum_{j=1}^pD_{X_j}(i_{N,j}T)(
X_1,\ldots,\hat{X}_j,\ldots,X_p)\\\di
+2\sum_{l<j}<X_l,X_j>T(X_1,\ldots,\stackrel{l}{\overbrace{N}},\ldots,
\stackrel{j}{\overbrace{N}},\ldots,X_p)
-N.N.T(X_1,\ldots,X_p)\\\di
\di+2\sum_{j=1}^pN.T(X_1,\ldots,D_{N}X_j,\ldots,X_p)
\di-\sum_{j=1}^pT(X_1,\ldots,D_{N}D_{N}X_j,\ldots,X_p)\\\di
-2\sum_{i<j}T(X_1,\ldots,D_{N}X_i,\ldots,D_{N}X_j,\ldots,X_p).\end{array}$$

Remark that, in restriction to $S^n$, the following equality holds
$$
\sum_{j=1}^pD_{X_j}(i_{N,j}T)( X_1,\ldots,\hat{X}_j,\ldots,X_p)=
\sum_{j=1}^pD_{X_j}(i_{\overrightarrow{r},j}T)(
X_1,\ldots,\hat{X}_j,\ldots,X_p).$$

 Now by using (7) and (11) and by taking the
restriction to $S^n$, we have
$$\begin{array}{l}\di
2\sum_{j=1}^pN.T(X_1,\ldots,D_{N}X_j,\ldots,X_p)=\\\di
2\sum_{j=1}^pN.T(X_1,\ldots,[{N},X_j],\ldots,X_p) +
2\sum_{j=1}^pN(\frac1r)T(X_1,\ldots,X_j,\ldots,X_p)\\\di
+2\sum_{j=1}^pN.T(X_1,\ldots,X_j,\ldots,X_p)
-2\sum_{j=1}^pN(<X_j,N>)T(X_1,\ldots,\stackrel{j}{\overbrace{N}},\ldots,X_p)=\\\di
2\sum_{j=1}^pN.T(X_1,\ldots,[{N},X_j],\ldots,X_p)
-2pT(X_1,\ldots,X_p)+2pN.T(X_1,\ldots,X_j,\ldots,X_p)\\\di
-2\sum_{j=1}^pN(<X_j,N>)T(X_1,\ldots,\stackrel{j}{\overbrace{N}},\ldots,X_p).\\\di
\sum_{j=1}^pT(X_1,\ldots,D_{N}D_{N}X_j,\ldots,X_p)=\\\di
\sum_{j=1}^pT(X_1,\ldots,[N,[N,X_j],\ldots,X_p)
-2\sum_{j=1}^pN(<X_j,N>)T(X_1,\ldots,\stackrel{j}{\overbrace{N}},\ldots,X_p).\end{array}$$
$$\begin{array}{l}\di
\sum_{i<j}T(X_1,\ldots,D_{N}X_i,\ldots,D_{N}X_j,\ldots,X_p)=\\\di
\sum_{i<j}T(X_1,\ldots,[{N},X_i],\ldots,[{N},X_j],\ldots,X_p)
+\frac{p(p-1)}2T(X_1,\ldots,X_p)\\\di
+\sum_{i<j}T(X_1,\ldots,X_i,\ldots,[{N},X_j],\ldots,X_p)+
\sum_{i<j}T(X_1,\ldots,[{N},X_i],\ldots,X_j,\ldots,X_p).\end{array}$$

So we get, in restriction to $S^n$
$$\begin{array}{l}\di
-N.N.T(X_1,\ldots,X_p)
\di+2\sum_{j=1}^pN.T(X_1,\ldots,D_{N}X_j,\ldots,X_p)\\
\di-\sum_{j=1}^pT(X_1,\ldots,D_{N}D_{N}X_j,\ldots,X_p)
-2\sum_{i<j}T(X_1,\ldots,D_{N}X_i,\ldots,D_{N}X_j,\ldots,X_p)=\\\di
-L_N\circ
L_NT(X_1,\ldots,X_p)+2pL_NT(X_1,\ldots,X_p)-p(1+p)T(X_1,\ldots,X_p).\end{array}$$

The curvature of $S^n$ is given by
$$
R(X,Y)Z=<X,Y>Z-<Y,Z>X\quad\mbox{and}\quad r(X)=(n-1)X.$$ Hence, a
direct computation gives that the curvature operator is given by
\begin{eqnarray*}
R(T)(X_1,\ldots,X_p)&=&p(n-1)T(X_1,\ldots,X_p)
+2\sum_{\sigma\in\T}T^\sigma(X_1,\ldots,X_p)\\
&&-2\sum_{i<j}\sum_{l=1}^n<X_i,X_j>T(X_1,\ldots,E_l,\ldots,E_l,\ldots,X_p).\end{eqnarray*}

Finally, we get
\begin{eqnarray*}
i^*(\De_{\reel^{n+1}}T)&=&\De_{S^n}i^*T+i^*\left(p(1-p)T+(2p-n)L_{N}T-
L_{N}\circ
L_{N}T\right.\\
&&\left.-2\sum_{\sigma\in\T}T^\sigma+O(T)\right),\end{eqnarray*}

One can conclude the proof by remarking that
$$i^*(L_NT)=i^*(L_{\overrightarrow{r}}T)\qquad\mbox{and}\quad
i^*(L_N\circ
L_NT)=-i^*(L_{\overrightarrow{r}}T)+i^*(L_{\overrightarrow{r}}\circ
L_{\overrightarrow{r}}T).$$ Q.E.D.

\begin{co} Let $\al$ be a differential  $p$-form on $\reel^{n+1}$. Then
\[
i^*(\De_{\reel^{n+1}}\al)
=\De_{S^n}i^*\al+i^*\left((2p-n+1)L_{\overrightarrow{r}}\al-L_{\overrightarrow{r}}\circ
L_{\overrightarrow{r}}\al-2di_{\overrightarrow{r}}\al\right)\]\end{co}

\begin{co} Let $T$ be a symmetric $p$-tensor field on $\reel^{n+1}$. Then
\begin{eqnarray*}
i^*(\De_{\reel^{n+1}}T)&=&\De_{S^n}i^*T+i^*\left(2p(1-p)T+(2p-n+1)L_{\overrightarrow{r}}T-
L_{\overrightarrow{r}}\circ
L_{\overrightarrow{r}}T\right.\nonumber\\
&&-\left.2{\de}^*(i_{\overrightarrow{r}}T)+2Tr(T)\odot<,>\right),\end{eqnarray*}where
$\odot$ is the symmetric product.
\end{co}

\section{ Eigenvalues and  eigenforms of $\De_{S^n}$ acting on
$\Om(S^n)$}

In this section, we will use corollary 2.1 and the results
developed in [10] to deduce  the eigenvalues and the spaces of
eigenforms of $\De_{S^n}$ acting on $\Om^*(S^n)$. We recover the
results of [10] and [11] in a more precise form.

 Let $\wedge^p H_k$ be the space of all coclosed harmonic homogeneous
$p$-forms of degree $k$ on $\reel^{n+1}$. A differential form
$\al$ belongs to $\wedge^pH_k$ if $\de(\al)=0$ and $\al$ can be
written
$$\al=\sum_{1\leq i_1<\ldots<i_p\leq n+1}\al_{i_1\ldots
i_p}dx_{i_1}\wedge\ldots\wedge dx_{i_p},$$where $\al_{i_1\ldots
i_p}$ are harmonic polynomial functions on $\reel^{n+1}$ of degree
$k$. For any $\al\in\wedge^pH_k$, we have
\begin{equation}
L_{\overrightarrow{r}}\al=di_{\overrightarrow{r}}\al+i_{\overrightarrow{r}}d\al=(k+p)\al.
\end{equation}

We have (see [10]),
$$i^*:\sum_{k\geq0} \wedge^pH_k\too\Om^p(S^n)$$is injective and its
image is dense.

For any $\al\in  \wedge^pH_k$,  we put
\begin{equation}
\om(\al)=\al-\frac1{p+k}di_{\overrightarrow{r}}\al.
\end{equation}

\begin{Le}

We get a linear map $\om:\wedge^pH_k\too \wedge^pH_k$ which is a
projector, i.e., $\om\circ\om=\om$. Moreover,
$$Ker\om=d(\wedge^{p-1}H_{k+1}),\quad
Im\om=\wedge^pH_k\cap Keri_{\overrightarrow{r}},$$and hence
$$\wedge^pH_k=\wedge^pH_k\cap Keri_{\overrightarrow{r}}\oplus
d(\wedge^{p-1}H_{k+1}).$$\end{Le}

The following lemma is an immediate consequence of Corollary 2.1
and (12).

 \begin{Le}\begin{enumerate}\item For any $\al\in  \wedge^pH_k\cap
Keri_{\overrightarrow{r}}$, we have
$$
\De_{S^n}i^*\al=(k+p)(k+n-p-1)i^*\al.$$ \item For any $\al\in
d(\wedge^{p-1}H_{k+1})$, we have
$$
\De_{S^n}i^*\al=(k+p)(k+n-p+1)i^*\al.$$
\end{enumerate}
\end{Le}
\begin{rem} We have
$$(k+p)(k+n-p-1)=(k'+p)(k'+n-p+1)\Leftrightarrow
k=k'+1\quad\mbox{and}\quad n=2p.$$\end{rem}

The following Table gives explicitly the spectra of $\De_{S^n}$
and the spaces of eigenforms with its  multiplicities . The
multiplicity was computed  in [11].\bigskip

\begin{center}{\bf Table I}\end{center}
 {\small
\begin{tabular}{|l|l|l|l|}
\hline $p$&The eigenvalues&The space of eigenforms&Multiplicity\\
\hline
&&&\\ $p=0$&$k(k+n-1),k\in\nat$&$\wedge^0 H_k$&$\frac{(n+k-2)!(n+2k-1)}{k!(n-1)!}$\\
&&&\\
\hline&&&\\
 $1\leq p\leq n,$
&$(k+p)(k+n-p-1),$ &
$\om(\wedge^pH_k)$&$\frac{(n+k-1)!(n+2k-1)}{p!(k-1)!(n-p-1)!(n+k-p-1)(k+p)}$\\
$n\not=2p$&$k\in\nat^*$&&\\ &$(k+p)(k+n-p+1),$
&$d(\wedge^{p-1}H_{k+1})$&$\frac{(n+k)!(n+2k+1)}{(p-1)!k!(n-p)!(n+k-p+1)(k+p)}$\\
&$k\in\nat$&&\\
\hline
&&&\\
 $1\leq p\leq n,$
&$(k+p)(k+p+1)$&
&\\
$n=2p$&$k\in\nat$&$\om(\wedge^p
H_{k+1})\oplus d(\wedge^{p-1}H_{k+1})$&$\frac{2(2p+k)!(2p+2k+1)}{p!(p-1)!k!(k+p+1)(k+p)}$\\
&&&\\
\hline

\end{tabular} }

\section{ Eigenvalues and  eigentensors of $\De_{S^n}$ acting on
$\s(S^n)$}

This section is devoted to the determination of the eigenvalues
and the spaces of  eigentensors of $\De_{S^n}$ acting on
$\s(S^n)$.

Let $\s^pP_k$ be the space of $T\in\s^p(\reel^{n+1})$ of the form
$$T=\sum_{1\leq i_1\leq\ldots\leq i_p\leq n+1}T_{i_1\ldots
i_p}dx_{i_1}\odot\ldots\odot dx_{i_p},$$where $T_{i_1\ldots i_p}$
are homogeneous polynomials of degree $k$. We put
$$
\s^pH_k^\de=\s^pP_k\cap Ker\De_{\reel^{n+1}}\cap
Ker\de\quad\mbox{and}\quad \s^pH_k^{\de0}=\s^pH_k^\de\cap Ker
Tr.$$

In a similar  manner as in [10] Lemma 6.4 and Corollary 6.6, we
have
\begin{equation}
\s^p P_k=\s^p H_k^\de\oplus(r^2\s^p P_{k-2}+dr^2\odot\s^{p-1}
P_{k-1}),\end{equation} and $$\di i^*:\sum_{k\geq0} \s^p
H_k^\de\too\s^pS^n$$ is injective and its image is dense in
$\s^pS^n$.

Now, for any $k\geq0$, we proceed to  give a direct sum
decomposition of $\s^pH_k^{\de}$ consisting of eigenspaces of
$\De_{S^n}$ and, hence, we determine completely the eigenvalues of
$\De_{S^n}$ acting on $\s^p(S^n)$. This will be done in several
steps.

 At  first, we have the following direct
sum decomposition:
\begin{equation}
\s^{p}H_k^\de=\s^pH_k^{\de0}\oplus\bigoplus_{l=1}^{[\frac{p}2]}
\s^{p-2l}H_k^{\de0}\odot<,>^l,
\end{equation}where $<,>^l$ is the symmetric product of $l$ copies
of $<,>$.

The task is now to decompose $\s^pH_k^{\de0}$ as a sum of
eigenspaces of $\De_{S^n}$ and get, according to (5), all the
eigenvalues. This decomposition  needs some preparation.

\begin{Le}Let $T\in\s^pP_k$ and $h\in\nat^*$. Then we have the following formulas:
\begin{enumerate}\item
$\de^*(i_{\overrightarrow{r}}T)-i_{\overrightarrow{r}}\de^*(T)=(p-k)T;$

\item$
\de^{*(h)}(i_{\overrightarrow{r}}T)-i_{\overrightarrow{r}}\de^{*(h)}(T)=h(p-k+h-1)
\de^{*(h-1)}(T);$
\item$\de^*(i_{\overrightarrow{r}^h}T)-i_{\overrightarrow{r}^h}\de^*(T)=
h(p-k-h+1)i_{\overrightarrow{r}^{h-1}}T,$\end{enumerate}where
$i_{\overrightarrow{r}^h}=\stackrel{h}{\overbrace{i_{\overrightarrow{r}}\circ\ldots
\circ i_{\overrightarrow{r}}}}$ and
$\de^{*(h)}=\stackrel{h}{\overbrace{\de^*\circ\ldots\circ\de^*}}.$\end{Le}

{\bf Proof.} The first formula is easily verified and the others
follow by induction on $h$. Q.E.D.\bigskip

Now, we will  construct two linear maps
$p_{\de^*}:\s^pP_k\too\s^pP_k$ for $k\leq p$, and
$p_{\overrightarrow{r}}:\s^pP_k\too\s^pP_k$ for $k\geq p$
satisfying:
\begin{enumerate}\item
$ p_{\de^*}\circ p_{\de^*}=p_{\de^*}$,
$Kerp_{\de^*}=i_{\overrightarrow{r}}(\s^{p+1}P_{k-1})$,
$Imp_{\de^*}=Ker\de^*\cap\s^pP_k$;

\item $ p_{\overrightarrow{r}}\circ
p_{\overrightarrow{r}}=p_{\overrightarrow{r}}$,
$Kerp_{\overrightarrow{r}}=\de^*(\s^{p-1}P_{k+1})$,
$Imp_{\overrightarrow{r}}=Keri_{\overrightarrow{r}}\cap\s^pP_k$.\end{enumerate}

The procedure is to    put, for $T\in \s^p P_k$,
\begin{eqnarray*}
p_{\de^*}(T)=\sum_{s=0}^k\al_si_{\overrightarrow{r}^s}\de^{*(s)}(T),\quad\mbox{and}\quad
p_{\overrightarrow{r}}(T)=\sum_{s=0}^p\be_s\de^{*(s)}(i_{\overrightarrow{r}^s}T),
\end{eqnarray*}
and find  $(\al_0,\ldots,\al_k)$ and $(\be_0,\ldots,\be_p)$ such
that the required properties are satisfied.

 A straightforward  computation using Lemma 4.1 gives
\begin{eqnarray*}
\de^*(p_{\de^*}(T))&=&\sum_{s=0}^{k-1}
(\al_s-(s+1)(k-p-s-2)\al_{s+1})i_{\overrightarrow{r}^{s}}\de^{*({s+1})}(T),\\
i_{\overrightarrow{r}}(p_{\overrightarrow{r}}(T))&=&\sum_{s=0}^{p-1}
(\be_s-(s+1)(p-k-s-2)\be_{s+1})\de^{*(s)}(i_{\overrightarrow{r}^{s+1}}T).
\end{eqnarray*}

Hence, we define $p_{\de^*}$ and $p_{\overrightarrow{r}}$ as
follows:
$$\left\{\begin{array}{lll}
p_{\de^*}(T)=&\di\sum_{s=0}^k\al_si_{\overrightarrow{r}^s}\de^{*(s)}(T)\\
\al_0=1&\mbox{and}\;
\al_s-(s+1)(k-p-s-2)\al_{s+1}=0&\;\mbox{for}\; 1\leq s\leq
k-1;\end{array}\right.$$

$$\left\{\begin{array}{lll}
p_{\overrightarrow{r}}(T)=&\di\sum_{s=0}^p\be_s\de^{*(s)}(i_{\overrightarrow{r}^s}T)\\
\be_0=1&\mbox{and}\;
\be_s-(s+1)(p-k-s-2)\be_{s+1}=0&\;\mbox{for}\; 1\leq s\leq
p-1.\end{array}\right.$$

From this definition and by using Lemma 4.1, one can check easily
 that
 $p_{\de^*}$ and $p_{\overrightarrow{r}}$ satisfy the required
 properties.

 On other hand, it is easy to check that we have, for any symmetric tensor
 field $T$ on $\reel^{n+1}$,
 \begin{eqnarray}
 \De_{\reel^{n+1}}(i_{\overrightarrow{r}}T)&=&i_{\overrightarrow{r}}\De_{\reel^{n+1}}(T)+
2\de T,\\
 \de(i_{\overrightarrow{r}}T)&=&i_{\overrightarrow{r}}\de(T)-Tr(T),\\
 Tr(\de^*(T))&=&-2\de(T)+\de^*(Tr(T)),\\
 Tr(i_{\overrightarrow{r}}T)&=&i_{\overrightarrow{r}}Tr(T).\end{eqnarray}From
 these
 formulas and from (3), one  deduce easily
 that $p_{\de^*}(\s^pH_k^{\de0})\subset
\s^pH_k^{\de0}$ and $p_{\overrightarrow{r}}(\s^pH_k^{\de0})\subset
\s^pH_k^{\de0}$ and thus one get the following direct sum
decompositions:
 \begin{eqnarray}
 \s^pH_k^{\de0}&=&\s^pH_k^{\de0}\cap Ker\de^*\oplus
 i_{\overrightarrow{r}}\left(\s^{p+1}H_{k-1}^{\de0}\right),\quad
 \mbox{if}\; k\leq p,\\
\s^pH_k^{\de0}&=&\s^pH_k^{\de0}\cap
Keri_{\overrightarrow{r}}\oplus
 \de^*\left(\s^{p-1}H_{k+1}^{\de0}\right),\quad
 \mbox{if}\; k\geq p.\end{eqnarray}

These decompositions are far for being sufficient and, in order to
obtain a more sharp direct sum decompositions of $\s^pH_k^{\de0}$,
we need the following lemma.

\begin{Le} \begin{enumerate}\item For $k<p$,
$i_{\overrightarrow{r}}:\s^pP_k\too\s^{p-1}P_{k+1}$ is injective.
\item For $k> p$, $\de^*:\s^pP_k\too\s^{p+1}P_{k-1}$ is injective.
\item For $k=p$,
$Ker\de^*=Ker{i_{\overrightarrow{r}}}.$\end{enumerate}\end{Le}

{\bf Proof.}\begin{enumerate}\item Let $T\in\s^pP_k$ such that
$i_{\overrightarrow{r}}T=0$. The second formula in  Lemma 4.1
gives, for any $h\geq1$,
$$i_{\overrightarrow{r}}\de^{*(h)}(T)=-h(p-k+h-1)
\de^{*(h-1)}(T).$$ Since $\de^{*(h)}(T)=0$ for $h\geq k+1$ and
$h(p-k+h-1)\not=0$ for any $h\geq1$, we get form this relation
that $\de^{*(h-1)}(T)=0$ for any $h\geq1$, in particular for
$h=1$, we get $T=0$. \item The same argument as 1. using the third
formula in Lemma 4.1. \item Let $T\in\s^pP_p$ such that
$i_{\overrightarrow{r}}T=0$. From Lemma 4.1, we get
$i_{\overrightarrow{r}}\de^*(T)=0$. Since
$\de^*(T)\in\s^{p+1}P_{p-1}$ and  from 1. we deduce that
$\de^*(T)=0$ and hence $Keri_{\overrightarrow{r}}\subset
Ker\de^*$. The same argument using Lemma 4.1 and 2. will give the
other inclusion. Q.E.D.\end{enumerate}

  By combining (20) and (21) with Lemma 4.2, we obtain  the following
  lemma.

\begin{Le} We have:
\begin{enumerate}\item if $k<p$
$$\s^pH_k^{\de0}=\bigoplus_{l=0}^k
i_{\overrightarrow{r}^{l}}\left(\s^{p+l}H_{k-l}^{\de0}\cap
Ker\de^*\right);$$ \item if $k>p$
$$
\s^pH_k^{\de0}=\bigoplus_{l=0}^{p}\de^{*l}\left(\s^{p-l}H_{k+l}^{\de0}\cap
Keri_{\overrightarrow{r}}\right);$$ \item If $k=p$,  for any
$0\leq l\leq p$,$$
i_{\overrightarrow{r}^l}\left(\s^{p+l}H_{p-l}^{\de0}\cap
Ker\de^*\right)=\de^{*l}\left(\s^{p-l}H_{p+l}^{\de0}\cap
Keri_{\overrightarrow{r}}\right),$$
and$$\s^{p}H_p^{\de0}=\bigoplus_{l=0}^p
i_{\overrightarrow{r}^{l}}\left(\s^{p+l}H_{p-l}^{\de0}\cap
Ker\de^*\right)=\bigoplus_{l=0}^{p}\de^{*l}\left(\s^{p-l}H_{p+l}^{\de0}\cap
Keri_{\overrightarrow{r}}\right).$$
\end{enumerate}
\end{Le}

Now, we use Corollary 2.2 to show that the decompositions of
$\s^pH_k^{\de0}$ given in Lemma 4.3 are composed by eigenspaces of
$\De_{S^n}$.

\begin{th}We have:\begin{enumerate}
 \item
If $k\leq p$, for any $0\leq q\leq k$ and any $T\in
i_{\overrightarrow{r}^{(k-q)}}\left(\s^{p+k-q}H_{q}^{\de0}\cap
Ker\de^*\right)$,
$$\De_{S^n}i^*T=\left((k+p)(n+p+k-2q-1)+2q(q-1)\right)i^*T;$$
\item If $k\geq p$, for any $0\leq q\leq p$ and for any
$T\in\de^{*(p-q)}\left(\s^{q} H_{k+p-q}^{\de0}\cap
Keri_{\overrightarrow{r}}\right)$,
$$\De_{S^n}i^*T=\left((k+p)(n+p+k-2q-1)+2q(q-1)\right)i^*T.
$$
\end{enumerate}\end{th}

{\bf Proof.}\begin{enumerate}\item Let
$T=i_{\overrightarrow{r}^{(k-q)}}(T_0)$ with
$T_0\in\s^{p+k-q}H_{q}^{\de0}\cap Ker\de^*$. We have from
Corollary 2.2
\begin{eqnarray*}
\De_{S^n}i^*T&=&i^*\left(2p(p-1)T+(n-2p-1)L_{\overrightarrow{r}}T+
L_{\overrightarrow{r}}\circ L_{\overrightarrow{r}}T\right.\\
&&\left.+2{\de}^*(i_{\overrightarrow{r}}T)-2Tr(T)\odot<,>\right).\end{eqnarray*}

We have
$$TrT=0,\quad L_{\overrightarrow{r}}=(k+p)T\quad\mbox{and}\quad
L_{\overrightarrow{r}}\circ L_{\overrightarrow{r}}T=(k+p)^2T.$$
Moreover, by using Lemma 4.1, we have
\begin{eqnarray*}
2{\de}^*(i_{\overrightarrow{r}}T)&=&2{\de}^*(i_{\overrightarrow{r}^{(k-q+1)}}T_0)\\
&\stackrel{\de^*(T_0)=0}=&2(k-q+1)(p+k-q-q-k+q-1+1)i_{\overrightarrow{r}^{(k-q)}}T_0\\
&=&2(k-q+1)(p-q)T.
\end{eqnarray*}
Hence
$$\De_{S^n}i^*T=(2p(p-1)+(n-2p-1)(k+p)+(k+p)^2+2(p-q)(k-q+1))i^*T.$$One
can deduce the desired relation by remarking  that
$$2p(p-1)+2(p-q)(k-q+1)=2(k+p)(p-q)+2q(q-1).$$
\item This follows by the same calculation as 1. Q.E.D.

\end{enumerate}

\bigskip

 From the fact that$$\di i^*:\sum_{k\geq0}
\s^p H_k^\de\too\s^pS^n$$ is injective and its image is dense in
$\s^pS^n$, from (15), and from  Lemma 4.3 and Theorem 4.1, note
that we have actually proved that
     the eigenvalues of $\De_{S^n}$ acting on $\s^pS^n$ belongs to
\begin{eqnarray*}
\left\{(k+p-2l)(n+p+k-2l-2q-1)+2q(q-1),\right.\\
\left.k\in\nat,0\leq l\leq[\frac{p}2],0\leq q\leq
min(k,p-2l)\right\}.\end{eqnarray*}Our next goal is to sharpen
this result by computing $\dim\s^{p}H_{k}^{\de0}\cap Ker\de^*$ if
$k\leq p$ and $\dim\s^{p}H_{k}^{\de0}\cap
Keri_{\overrightarrow{r}}$ if $k\geq p$.

\begin{Le} We have the following formulas:
\begin{enumerate}\item
$\dim\s^p H_k^\de=\dim\s^p P_{k}-\dim\s^p P_{k-2}-\dim\s^{p-1}
P_{k-1}+\dim\s^{p-1} P_{k-3},$ \item$\dim\s^{p}H_k^{\de0}=\dim\s^p
H_k^\de-\dim\s^{p-2} H_k^\de,$\item$ \dim(\s^{p}H_k^{\de0}\cap
Ker\de^*)=\dim\s^{p}H_k^{\de0}-\dim\s^{p+1}H_{k-1}^{\de0}\; (k\leq
p),$\item$ \dim(\s^{p}H_k^{\de0}\cap
Keri_{\overrightarrow{r}})=\dim\s^{p}H_k^{\de0}-\dim\s^{p-1}H_{k+1}^{\de0}\;(k\geq
p).$
\end{enumerate}
Note that we use the convention that $\s^p P_{k}=\s^p
H_k^\de=\s^{p}H_k^{\de0}=0$ if $k<0$ or $p<0$.
\end{Le}

{\bf Proof.} \begin{enumerate}\item The formula is a consequence
of (14), the relation $$(r^2\s^pP_{k-2})\cap
(dr^2\odot\s^{p-1}P_{k-1})=r^2(dr^2\odot\s^{p-1}P_{k-3})$$ and the
fact that $dr^2\odot.:\s^{p}P_{k}\too\s^{p+1}P_{k+1}$ is
injective. \item The formula is a consequence of (15). \item The
formula is a consequence of (20) and Lemma 4.2. \item The formula
is a consequence of (21) and Lemma 4.2. Q.E.D.
\end{enumerate}
\bigskip

A straightforward calculation using  Lemma 4.4 and the formula
$$\dim\s^p P_{k}=\frac{(n+p)!}{n!p!}\frac{(n+k)!}{n!k!}$$gives
$\dim\s^{p}H_{k}^{\de0}\cap Ker\de^*$ if $k\leq p$ and
$\dim\s^{p}H_{k}^{\de0}\cap Keri_{\overrightarrow{r}}$ if $k\geq
p$. We summarize the results on the following Table.\eject

\begin{center} {\bf Table II}\end{center}
\begin{tabular}{|c|c|c|}
\hline Space&Dimension& Conditions on $k$ and $p$\\
\hline
 $\s^0H_k^{\de0}\cap
Keri_{\overrightarrow{r}}$&$\di\frac{(n+k-2)!(n+2k-1)}{k!(n-1)!}$&$k\geq0$\\
\hline
$\s^pH_0^{\de0}\cap Ker\de^*$&$\di\frac{(n+p-2)!(n+2p-1)}{p!(n-1)!}$& $p\geq0$\\
\hline
 $\s^1H_k^{\de0}\cap
Keri_{\overrightarrow{r}}$&$\di\frac{(n+k-3)!k(n+2k-1)(n+k-1)}{
(n-2)!(k+1)!}$&
$k\geq1$\\
\hline $\s^pH_1^{\de0}\cap
Ker\de^*$&$\di\frac{(n+p-3)!p(n+2p-1)(n+p-1)}{
(n-2)!(p+1)!}$&$p\geq1$\\
\hline
 $\s^pH_k^{\de0}\cap
Ker\de^*$
&$\di\frac{(n+k-4)!(n+p-3)!(n+p+k-2)}{k!(p+1)!(n-1)!(n-2)!}\times$&\\
&$(n-2)(n+2k-3)(n+2p-1)(p-k+1)$&$2\leq k\leq p$\\
\hline $\s^pH_k^{\de0}\cap Keri_{\overrightarrow{r}}$
&$\di\frac{(n+k-3)!(n+p-4)!(n+p+k-2)}{(k+1)!p!(n-1)!(n-2)!}\times$&\\
&$(n-2)(n+2k-1)(n+2p-3)(k-p+1)$& $k\geq p\geq 2$\\
\hline
\end{tabular}

\begin{rem} Note that, for $n=2$, we have
\begin{eqnarray*}
\dim(\s^pH_k^{\de0}\cap Ker\de^*)=0\quad \mbox{for}\quad 2\leq
k\leq p,\\
\dim(\s^pH_k^{\de0}\cap Keri_{\overrightarrow{r}})=0\quad
\mbox{for}\quad k\geq p\geq 2.\end{eqnarray*}\end{rem}

For simplicity we introduce the following notations.
\begin{eqnarray*}
S_0&=&\left\{(k,l,q)\in\nat^3,0\leq l\leq[\frac{p}2],0\leq k\leq
p-2l,0\leq q\leq k\right\},\\
S_1&=&\left\{(k,l,q)\in\nat^3,0\leq l\leq[\frac{p}2], k>
p-2l,0\leq q\leq p-2l\right\},\\
V_{q,l}^k&=&i_{\overrightarrow{r}^{k-q}}\left(\s^{p-2l+k-q}H^{\de0}_q\cap
Ker\de^*\right)\odot<,>^l\;\quad\mbox{for} \quad(k,l,q)\in S_0,\\
W_{q,l}^k&=&\de^{*{(p-2l-q)}}\left(\s^qH^{\de0}_{p-2l+k-q}\cap
Keri_{\overrightarrow{r}}\right)\odot<,>^l\;\quad\mbox{for}\quad
(k,l,q)\in S_1.\end{eqnarray*}

Let us summarize all the results above.

\begin{th}\begin{enumerate}\item For $n=2$, we have:
\begin{enumerate} \item The set of the eigenvalues of $\De_{S^2}$
acting on $\s^pS^2$ is
$$\left\{(k+p-2l)(p+k-2l+1),\;k\in\nat,
0\leq l\leq[\frac{p}2]\right\};$$ \item The eigenspace associated
to the eigenvalue $\la(k,l)=(k+p-2l)(k+p-2l+1)$ is given by
$$V_{\la(k,l)}=\left\{\begin{array}{lll}
\di\bigoplus_{a=0}^{\min(l,[\frac{k}2])}\left(V_{0,l-a}^{k-2a}
\oplus
V_{1,l-a}^{k+1-2a}\right)&\mbox{if}&0\leq k\leq p-2l\\
\di\bigoplus_{a=0}^{\min(l,[\frac{k}2])}\left(W_{0,l-a}^{k-2a}
\oplus W_{1,l-a}^{k+1-2a}\right)&\mbox{if}&k>p-2l;
\end{array}\right.$$
\item The multiplicity of $\la(k,l)$ is given by
$$m(\la(k,l))=2(\min(l,[\frac{k}2])+1)(1+2p+2k-4l).$$

\end{enumerate}
\item For $n\geq3$, we have:\begin{enumerate}\item The set of the
eigenvalues of $\De_{S^n}$ acting on $\s^pS^n$ is
\begin{eqnarray*}
\left\{(k+p-2l)(n+p+k-2l-2q-1)+2q(q-1),\right.\\
\left.k\in\nat,0\leq l\leq[\frac{p}2],0\leq q\leq
min(k,p-2l)\right\};\end{eqnarray*} \item The space
$$\p=\sum_{k\geq0}\s^pH_k^{\de}=(\bigoplus_{(k,l,q)\in
S_0}V_{q,l}^k)\oplus(\bigoplus_{(k,l,q)\in S_1}W_{q,l}^k)$$ is
dense in $\s^pS^n$ and, for any $(k,q,l)\in S_0$ (resp.
$(k,q,l)\in S_1$), $V_{q,l}^k$ (resp. $W_{q,l}^k$) is a subspace
of the eigenspace associated to the eigenvalue
$(k+p-2l)(n+p+k-2l-2q-1)+2q(q-1)$; \item The dimensions of
$V_{q,l}^k$ and $W_{q,l}^k$ are given in Table II since
\begin{eqnarray*}\dim V_{q,l}^k&=&\dim\left(\s^{p-2l+k-q}H^{\de0}_q\cap
Ker\de^*\right)\;\quad\mbox{for} \quad(k,l,q)\in S_0,\\
\dim W_{q,l}^k&=&\dim\left(\s^qH^{\de0}_{p-2l+k-q}\cap
Keri_{\overrightarrow{r}}\right)\;\quad\mbox{for}\quad (k,l,q)\in
S_1.\end{eqnarray*}
\end{enumerate}\end{enumerate}\end{th}

{\bf References}\bigskip

[1] {\bf E. Bedford and T. Suwa,} Eigenvalues of Hopf manifolds,
American Mathemaical Society, Vol. {\bf 60} (1976), 259-264.

[2] {\bf B. L. Beers and R. S. Millman,} The spectra of the
Laplace-Beltrami operator on compact, semisimple Lie groups, Amer.
J. Math., {\bf 99} (4) (1975), 801-807.

[3] {\bf M. Berger and D. Ebin,} Some decompositions of the space
of symmetric tensors on Riemannian manifolds, J. Diff. Geom., {\bf
3} (1969), 379-392.

[4] {\bf M. Berger, P. Gauduchon and E. Mazet,} Le spectre d'une
vari\'et\'e riemannienne, Lecture Notes in Math., Vol {\bf 194},
Springer Verlag (1971).

[5] {\bf A. Besse,} Einstein manifolds, Springer-Verlag,
Berlin-Hiedelberg-New York (1987).

[6] {\bf  M. Boucetta ,} Spectre des Laplaciens de Lichnerowicz
sur les sph\`eres et les projectifs r\'eels, Publicacions
Matem\`atiques, Vol. {\bf 43} (1999), 451-483.

[7] {\bf  M. Boucetta ,} Spectre du Laplacien de Lichnerowicz sur
les projectifs complexes, C. R. Acad. Sci. Paris, t. 333, S\'erie
I, (2001), 571-576.

[8] {\bf S. Gallot and D. Meyer,} Op\'erateur de courbure et
laplacien des formes diff\'erentielles d'une vari\'et\'e
riemannienne, J. Math. Pures Appl., {\bf 54} (1975), 259-289.

[9] {\bf G. W. Gibbons and M. J. Perry,} Quantizing gravitational
instantons, Nuclear Physics B, Vol. {\bf 146}, Issue I (1978),
90-108.

[10] {\bf A. Ikeda and Y. Taniguchi,} Spectra and eigenforms of
the Laplacian on $S^n$ and $P^n(\comp)$, Osaka J. Math., {\bf 15}
(3) (1978), 515-546.

[11] {\bf I. Iwasaki and K. Katase,}  On the spectra of Laplace
operator on $\wedge^*(S^n)$, Proc. Japan Acad., {\bf 55}, Ser. A
(1979), 141-145.

[12] {\bf E. Kaneda,} The spectra of 1-forms on simply connected
compact irreducible Riemannian symmetric spaces, J. Math. Kyoto
Univ., {\bf 23} (1983), 369-395 and {\bf 24} (1984), 141-162.

[13] {\bf A. L\'evy-Bruhl-Laperri\`ere,} Spectre de de Rham-Hodge
sur les formes de degr\'e 1 des sph\`eres de $\reel^n$ ($n\geq6$),
Bull. Sc. Math., $2^e$ s\'erie, {\bf 99} (1975), 213-240.

[14] {\bf A. L\'evy-Bruhl-Laperri\`ere,} Spectre de de Rham-Hodge
sur l'espace projectif  complexe, C. R. Acad. Sc. Paris {\bf 284}
S\'erie A (1977), 1265-1267.

[15] {\bf A. Lichnerowicz,} Propagateurs et commutateurs en
relativit\'e g\'en\'erale, Inst. Hautes Etude Sci. Publ. Math.,
{\bf 10} (1961).

[16] {\bf K. Mashimo,} Spectra of Laplacian on $G_2/SO(4)$, Bull.
Fac. Gen. Ed. Tokyo Univ. of Agr. and Tech. {\bf 26} (1989),
85-92.

[17] {\bf K. Mashimo,} On branching theorem of the pair
$(G_2,SU(3))$, Nihonkai Math. J., Vol. {\bf 8} No. 2 (1997),
101-107.

[18] {\bf K. Mashimo,} Spectra of the Laplacian on the Cayley
projective plane, Tsukuba J. Math., Vol. {\bf 21} No. 2 (1997),
367-396.

[19] {\bf R. Michel,} Probl\`eme d'analyse g\'eom\'etrique li\'es
\`a la conjecture de Blaschke, Bull. Soc. Math. France, {\bf 101}
(1973), 17-69.

[20] {\bf K. Pilch and N. Schellekens,} Formulas of the
eigenvalues of the Laplacian on tensor harmonics on symmetric
coset spaces, J. Math. Phys., {\bf 25} (12) (1984), 3455-3459.

[21] {\bf C. Tsukamoto, } The sepctra of the Laplace-Beltrami
operator on $SO(n+2)/SO(2)\times SO(n)$ and $Sp(n+1)/Sp(1)\times
Sp(n)$, Osaka J. Math. {\bf 18} (1981), 407-226.

[22] {\bf N. P. Warner,} The spectra of operators on $\comp P^n$,
Proc. R. Soc. Lond. A {\bf 383} (1982), 217-230.

\end{document}